\documentclass[12pt]{article}
\usepackage{latexsym}
\usepackage{amsfonts}
\def\R{\mathbf{R}}
\def\N{\mathbf{N}}
\def\Ima{\mathrm{Im}\ }
\def\dist{\mathrm{dist}}
\author{Alexandre Eremenko\thanks{Supported by NSF grant
DMS-0555279.}}
\title{A version of Fabry's theorem for power series
with regularly varying coefficients}

\begin{document}
\maketitle

\begin{abstract} For real power series whose non-zero
coefficients satisfy $|a_m|^{1/m}\to~1$,
we prove a stronger version of Fabry's theorem
relating the frequency of sign changes in the
coefficients and analytic continuation of the
sum of the power series.
\end{abstract}

For a set $\Lambda$ of non-negative integers, we consider
the counting function
$$n(x,\Lambda)=\#\Lambda\cap[0,x].$$
We say that $\Lambda$ is {\em measurable} if the
limit
$$\lim_{x\to+\infty} n(x,\Lambda)/x$$
exists, and call this limit the {\em density} of $\Lambda$.

Let $S=\{ a_m\}$ be a sequence of real numbers. We say that
{\em a sign change} occurs at the place $m$
if there exists $k<m$ such that $a_ma_k<0$ while
$a_j=0$ for $k<j<m$. 
\vspace{.1in}

\noindent
{\bf Theorem A.} {\em The following two properties
of a set $\Lambda$ of positive integers are equivalent:
\vspace{.07in}

\noindent
(i) Every power series
\begin{equation}
\label{ps}
f(z)=\sum_{m=0}^\infty a_mz^m
\end{equation}
of radius of convergence $1$, with real coefficients
and such that the changes
of sign of $\{ a_m\}$ occur only
for $m\in\Lambda$, has a singularity
on the arc $$I_\Delta=\{ e^{i\theta}:|\theta|\leq\Delta\},$$
and
\vspace{.1in}

\noindent
(ii) For every $\Delta'>\Delta$ there exists
a measurable set $\Lambda'\subset\N$ of density $\Delta'$
such that $\Lambda\subset\Lambda'$.}
\vspace{.1in}

Implication (ii) $\longrightarrow$ (i) is a consequence
of Fabry's General Theorem \cite{Fabry,Bieber},
as restated by P\'olya.
For the implication (i) $\longrightarrow$ (ii)
see \cite{Koosis}. Fabry's General theorem takes into
account not only the sign changes of coefficients
but also the absolute values of coefficients.
It has a rather complicated statement and the sufficient
condition of the existence of a singularity given by
this theorem is not the best possible.
The best possible condition in Fabry's General theorem
is unknown, see, for example the discussion in \cite{E}.

Alan Sokal (private communication) asked
what happens if we assume that the
power series (\ref{ps})
satisfies the additional regularity condition:
\begin{equation}\label{0}
\lim_{m\in P, m\to\infty}|a_m|^{1/m}=1,
\end{equation}
where $P=\{ m:a_m\neq 0\}$.
This condition holds for most
interesting generating functions.
The answer is somewhat surprising:
\vspace{.1in}

\noindent
{\bf Theorem 1.} {\em The following two properties of
a set $\Lambda$ of positive integers are equivalent:
\vspace{.07in}

\noindent
a) Every power series $(\ref{ps})$ satisfying $(\ref{0})$,
with real coefficients
and such that the changes
of sign of the coefficients $a_m$ occur only
for $m\in\Lambda$, has a singularity
on the arc $I_\Delta$,
and
\vspace{.1in}

\noindent
b) 
All measurable subsets $\Lambda'\subset\Lambda$
have densities at most $\Delta$.}
\vspace{.1in}

We recall 
that the {\em minimum density}
$$\underline{D}_2(\Lambda)=\lim_{r\to0+}
\liminf_{x\to+\infty}\frac{n((r+1)x,\Lambda)-n(x,\Lambda)}{rx}$$
can be alternatively defined as
the $\sup$ of the limits
\begin{equation}\label{lim}
\lim_{x\to\infty}n(x,\Lambda')/x
\end{equation}
over all measurable sets $\Lambda'\subset\Lambda$.

Similarly the {\em maximum density} of $\Lambda$
is
$$\overline{D}_2(\Lambda)=\lim_{r\to 0+}\limsup_{x\to\infty}
\frac{n((r+1)x,\Lambda)-n(x,\Lambda)}{rx},$$
and it equals to inf of the limits (\ref{lim}) over all
measurable sequences of non-negative integers $\Lambda'$
containing $\Lambda$.

For all these properties of minimum and maximum
densities see \cite{P}.

Thus condition (ii) is equivalent to $\overline{D}_2(\Lambda)\leq\Delta$ while condition b) is equivalent to
$\underline{D}_2(\Lambda)\leq\Delta.$
\vspace{.1in}

\noindent
{\bf Corollary 1.} {\em The following
two properties of a set $\Lambda$ of positive
integers are equivalent:
\vspace{.07in}

\noindent
$A$. Every power series 
\begin{equation}
\label{pps}
\sum_{m\in\Lambda}a_mz^m
\end{equation}
satisfying $(\ref{0})$ has a singularity on
$I_\Delta$, 
\vspace{.07in}

\noindent
$A'$. Every power series $(\ref{pps})$ 
satisfying
$(\ref{0})$ has a singularity on every closed arc
of length $2\pi\Delta$ of the unit circle, and
\vspace{.07in}

\noindent 
$B.$ $\underline{D}_2(\Lambda)\leq\Delta.$}
\vspace{.1in}

Indeed, all assumptions
of $A$ are invariant with respect to
the change of
the variable $z\mapsto\lambda z,\; |\lambda|=1$, thus $A$
is equivalent to the formally stronger statement $A'$.

Now, the number of sign changes of any sequence
does not exceed the number of its non-zero terms, thus
$B$ imlies $A$ by Theorem~1. The remaining implication
$A\longrightarrow B$ will be proved in the end
of the proof of Theorem~1.
\vspace{.1in}

{\em Proof of Theorem 1.}
b) $\longrightarrow$ a). Proving this by contradiction,
we assume that $\underline{D}_2(\Lambda)\leq\Delta$, and there exists
a function $f$ of the form (\ref{ps})
with the property (\ref{0}) which has
an analytic continuation to $I_\Delta$,
and such that the sign changes occur only for $m\in\Lambda$.

Without loss of generality
we assume that $a_0=1$, and $\Delta<1$.
\vspace{.1in}

\noindent
{\bf Lemma 1.} 
{\em For a function $f$ as in $(\ref{ps})$ to
have an immediate analytic continuation from the unit disc
to the arc
$I_\Delta$ it is necessary and sufficient
that there exists an entire function $F$ of exponential type
with the properties
\begin{equation}
\label{3}
a_m=(-1)^mF(m),\quad\mbox{for all}\quad m\geq 0,
\end{equation}
and
\begin{equation}
\label{4}
\limsup_{t\to\infty}\frac{\log|F(te^{i\theta})|}{t}
\leq \pi b|\sin\theta|,\quad |\theta|<\alpha,
\end{equation}
with some $b<1-\Delta.$}
\vspace{.1in}

This result can be found in \cite{Arak}, see also 
\cite{Ar,E}.

Consider the sequence of subharmonic functions
$$u_m(z)=\frac{1}{m}\log|F(mz)|,\quad m=1,2,3,\ldots.$$
This sequence is uniformly bounded from above on
every compact subset of the plane, because $F$
is of exponential type. Moreover, 
$u_m(0)=0$ because of our assumption that $a_0=F(0)=1$.
Compactness Principle \cite[Th. 4.1.9]{Hor} implies that
from every sequence of integers $m$ one can choose a
subsequence such that the limit $u=\lim u_m$
exists. This limit is a subharmonic function in the plane
that satisfies in view of (\ref{4})
\begin{equation}
\label{5}
u(re^{i\theta})\leq\pi br|\sin\theta|,\quad |\theta|<\alpha,
\end{equation}
with some $b$ satisfying
$$0<b<1-\Delta.$$

We use the following result of P\'olya
\cite[footnote 18, p. 703]{P2}:
\vspace{.1in}

\noindent
{\bf Lemma 2.} {\em Let $f$ be a power series
$(\ref{ps})$ of radius of convergence $1$.
Let $\{ a_{m_k}\}$ be a subsequence of coefficients
with the property
$$\lim_{k\to\infty}|a_{m_k}|^{1/m_k}= 1,$$
and assume that for some $r>0$
the number of non-zero coefficients
$a_j$ on the interval
$m_k\leq j\leq(1+r)m_k$ is $o(m_kr)$ as $k\to\infty$.
Then $f$ has no analytic continuation to any point
of the unit circle.} 
\vspace{.1in}

Lemma 2 also follows from the results of \cite{Arak}
or \cite{E}.

Now we show that (\ref{0}) implies the following:
\vspace{.1in}

\noindent
{\bf Lemma 3.} {\em Every limit function has the property
$u(x)=0$ for $x\geq 0$.}
\vspace{.1in}

{\em Proof of Lemma 3.}
Let $U=\{ x:x\geq 0, u(x)<0\}$. This set is open
because $u$ is upper semi-continuous. 
Take any closed interval $J=[c,d]\subset U$.
Then
$u(x)\leq-\epsilon,\; x\in J,$ with some $\epsilon>0$.
Let $\{ m_k\}$ be the sequence of integers
such that $u_{m_k}\to u$. 
Then from the definition of $u_m$
we see that
$$\log|F(m_kx)|\leq -m_k\epsilon/2\quad\mbox{for}\quad
x\in J$$
and for all large $k$. Together with (\ref{3}) and (\ref{0})
this implies that $a_j=0$ for all 
$j\in m_kJ$.
Let $a_{m_k^\prime}$ be the last non-zero coefficient
before $cm_k$.
Applying Lemma 2 to the sequence $\{ m_k^\prime\}$
we conclude that $f$ has no analytic continuation
from the unit disc. This is a contradiction
which proves Lemma~3.
\hfill$\Box$
\vspace{.1in}

Now we use the following general fact:
\vspace{.1in}

\noindent
{\bf Grishin's Lemma.} {\em Let $u\leq v$ be two subharmonic
functions, and $\mu$ and $\nu$ their respective Riesz
measures. Let $E$ be a Borel set such that
$u(z)=v(z)>-\infty$ for $z\in E$.
Then the restrictions of the Riesz measures on $E$
satisfy
$$\mu\vert_E\leq \nu\vert_E.$$
}

The references are \cite{VP,G,F}.

In view of Lemma 2, we can apply
Grishin's Lemma to $u$ and $v(z)=\pi b|\Ima z|$
and $E=[0,\infty)\subset\R$. We obtain that
the Riesz measure $d\mu$ of any limit function $u$
of the sequence $\{ u_k\}$ satisfies
\begin{equation}
\label{dmu}
d\mu\vert_{[0,\infty)}\leq b\ dx.
\end{equation}
Now we go back to our coefficients and function $F$.
By our assumption, the sign changes
occur on a sequence $\Lambda$ whose minimum density
is at most $\Delta$.
Choose a number $a$ such that
$b<a<1-\Delta$.
By the first definition of the minimum density,
there exist $r>0$ and a sequence $x_k\to\infty$
such that 
$$n((1+r)x_k,\Lambda)-n(x_k,\Lambda)\leq (1-a)rx_k.$$

\noindent
{\bf Lemma 4.}
{\em Let $(a_0,a_1,\ldots,a_N)$ be a sequence of
real numbers, and $f$ a real analytic function on the closed interval
$[0,N]$, such that $f(n)=(-1)^na_n$. Then the number of zeros
of $f$ on $[0,N]$, counting multiplicities, is at least
$N$ minus the number of sign changes of the sequence $\{ a_n\}$.}
\vspace{.1in}

{\em Proof.} Consider first an interval $(k,n)$ such that
$a_ka_n\neq 0$ but $a_j=0$ for $k<j<n$.
We claim that $f$ has at least
$$n-k-\#(\mbox{sign changes in the pair}\; (a_k,a_n))$$
zeros
on the open interval $(k,n)$. Indeed, the number of zeros of $f$
on this interval is at least $n-k-1$ in any case.
This proves the claim if there is a sign change in the pair $(a_k,a_n)$.
If there is no sign change, that is $a_na_k>0$, then  $f(n)f(k)=(-1)^{n-k}$.
So the number of zeros of $f$ on the interval $(n,k)$ is of the same
parity as $n-k$. But $f$ has at least $n-k-1$ zeros on this
interval, thus the total number of zeros is at least $n-k$.
This proves our claim.

Now let $a_k$ be the first and $a_n$ the last non-zero term of our sequence.
As the interval $(k,n)$ is a disjoint union of the intervals
to which the above claim applies, we conclude that the number of
zeros of $f$ on $(k,n)$ is at least $(n-k)$ minus the number of
sign changes of our sequence. On the rest of the interval
$[0,N]$
our function has at least $N-n+k$ zeros, so the total
number of zeros is at least $N$ minus the number of sign changes.
\hfill$\Box$
\vspace{.1in}

Let $u$ be a limit function of the subsequence
$\{ u_{m_k}\}$ with $m_k=[x_k]$.
By Lemma 4, the function $F$ has at least
$arx_k-2$ zeros on each interval $[x_k,(1+r)x_k]$, 
which implies that the Riesz measure $\mu$ of $u$
satisfies
$$\mu([1,1+r])\geq ar.$$
This contradicts (\ref{dmu}) and thus proves 
the implication b) $\longrightarrow$ a).
\vspace{.1in}

a) $\longrightarrow$ b). Suppose that a set
$\Lambda$ of positive integers does not
satisfy b).
We will construct power series $f$ of the form (\ref{pps})
which has an immediate analytoc continuation from
the unit disc to the arc $I_\Delta$. This will simultaneously
prove the implications a) $\longrightarrow$ b) of
Theorem~1 and $A\longrightarrow B$ of Corollary~1.
 
Let $\Lambda'\subset\Lambda$ be a measurable set
of density $\Delta'>\Delta$. Let $S$
be the complement of
$\Lambda'$ in the set of positive integers. Then
$S$ is also measurable and has density $1-\Delta'$.

Consider the infinite product
$$F(z)=\prod_{t\in S}\left(1-\frac{z^2}{t^2}\right).$$
This is an entire function of exponential type
with indicator $\pi(1-\Delta')|\sin\theta|$, and furthermore,
\begin{equation}
\label{last}
\log|F(z)|\geq \pi(1-\Delta')|\Ima z|+o(|z|),
\end{equation}
as $z\to\infty$
outside the set $\{ z:\dist (z,S)\leq 1/4\}$.
(See \cite[Ch. II, Thm. 5]{L} for this result.)
Now we use the sufficiency part of Lemma~1,
and define the coefficients of our power series
by $a_m=(-1)^mF(m)$. Then we have all needed properties,
in particular (\ref{0}) follows from (\ref{last}).
\vspace{.1in}

The author thanks Alan Sokal for many illuminating conversations
about Fabry's theorem.

\vspace{.1in}

{\em Purdue University,

West Lafayette, Indiana, USA

eremenko@math.purdue.edu}

\begin{thebibliography}{11}
\bibitem{Arak} 
 N. U. Arakelyan and V. A. Martirosyan,
Localization of singularities on the boundary of the
circle of convergence,
Izvestiya Akademii Nauk Armyanskoi SSR,
Mat. vol. 22, No. 1 (1987) 3-21 (Russian)
English translation:
Journal of Contemporary Mathematical Analysis,
22 (1988) 1--19.
\bibitem{Ar} N. Arakelyan, W. Luh and J. M\"uller,
On the localization of singularities of lacunar
power series, Complex Variables and Elliptic Equations,
52 (2007) 651--573.
\bibitem{Bieber} L. Bieberbach,
 Analytische Fortsetzung, Springer,
Berlin, 1955.
\bibitem{E} A. Eremenko, Densities in Fabry's theorem,
preprint arXiv:0709.2360.
\bibitem{F} B. Fuglede,
Some properties of the Riesz charge associated with
a $\delta$-subharmonic function.  
Potential Anal. 1 (1992) 355--371. 
\bibitem{Fabry} E. Fabry,
Sur les s\'eries de Taylor qui ont
une infinit\'e de points singuliers,
Acta math., 22 (1898) 65--87.
\bibitem{G} A. F. Grishin,
Sets of regular growth of entire functions I,  
Teor. Funktsii Funktsional. Anal. i Prilozhen.
No. 40 (1983) 36--47 (Russian). 
\bibitem{Hor} L. H\"ormander, The analysis of linear partial
differential operators, vol. I, Springer, Berlin, 1983.
\bibitem{Koosis} P. Koosis, The Logarithmic Integral, vol. II
Cambridge Univ. Press, Cambridge, 1992.
\bibitem{L} B. Ya. Levin, Distribution of zeros of entire
functions, AMS Providence, RI, 1980.
\bibitem{P2} G. P\'olya, \"Uber gewisse notwendige
Determinantkriterien fur die
Fortsetzbarkeit einer Potenzreiche,
Math. Ann. 99 (1928) 687--706.
\bibitem{P} G. P\'olya, Untersuchungen \"uber L\"ucken und
Singularit\"aten von Potenzreichen, Math. Z.,
29 (1929) 549--640.
\bibitem{VP} C. de la Valle-Poussin, Potentiel et
probl\'eme g\'en\'eralis\'e de Dirichlet,
Math. Gazette, 22 (1938) 17--36.
\end{thebibliography}
\end{document}